\documentclass{article}
\usepackage{amsmath}
\usepackage{amssymb}
\usepackage{amsthm}
\newtheorem{lemma}{Lemma}
\newtheorem{definition}{Definition}
\newtheorem{theorem}{Theorem}
\author{L.C. Brown}
\title{Definability of the Integrability Locus in Polynomially Bounded o-Minimal Structures}
\begin{document}
\maketitle
	Let $\tilde{R}$ be an o-minimal expansion of the real number field $\mathbf{R}$. Assume also that $\tilde{R}$ is polynomially bounded: that is, if $f:\mathbf{R}\rightarrow\mathbf{R}$ is definable in $\tilde{R}$, there is some integer $n$ and some real number $C$ such that, whenever $x > C$, $f(x) < x^n.$ By the cell decomposition theorem in \cite{pill}, every function definable in $\tilde{R}$ is measurable. If $f:\mathbb{R}^{m+n}\rightarrow\mathbb{R}$ is a non-negative function definable in $\tilde{R}$, we may define its integrability locus $$L(f) = \{x\in\mathbb{R}^m: \int_{\mathbb{R}^n} f(x,y) dy <\infty\}.$$
\\
We claim $L(f)$ is definable in $\tilde{R}.$ For ease of notation, we instead consider, for $A\subset\mathbb R^{m+n}$, the set $\{x\in\mathbb R^m: \mathrm{vol}(A_x)\}<\infty$, where $A_x = \{y\in\mathbb R^n: (x,y)\in A\}.$ Indeed, we do so by piecing together pieces already in the literature. We need the strong parameter independence of the definition of Lipschitz stratifications from Halupczok and Yin in \cite{halup}, the Lipschitz stratification in polynomially bounded o-minimal structures obtained by Nguyen in \cite{nguyen}, the Lipschitz Thom-Mather lemma obtained by Parusinski in \cite{parus}, and the volume estimates on inversions about the unit sphere appearing in Comte, Lion, and Rolin in \cite{clr}. We need to fix the proof of a lemma in Halupczok and Yin, clarify the statement of Nguyen's result, and explicate Parusinski's argument. Some experts hold that the result we prove here is known; however, the complete proof is not written in full anywhere, and the missing components were not obvious to the author.
\section{Valuative Lipschitz Stratifications}
Halupczok and Yin examine o-minimal expansions of fields containing $\mathbb{R}$ in \cite{halup}. We only need their Proposition 1.2.5. The non-trivial part of this proof relies on Proposition 1.6.11, namely the implication from (2') to (3') in its statement. We adopt their notation and definitions.
\\ \\
Let $\mathcal{R}$ be a polynomially bounded o-minimal expansion of a real-closed field. We write $d(a,X)$ for the distance between a point $a\in\mathcal R^m$ and a set $X\subset \mathcal R^m$; $d(a,\emptyset) = 1.$
\begin{definition}[\cite{halup},Definition 1.2.1; \cite{nguyen}, \S 4.2]
	Let $X\subset\mathcal R^m$ be definable with parameters and fix a definable $C^2$ cell decomposition $\mathcal C$ of $X.$ Write $\mathcal X = (X^i)_{i=0}^m,$ where $X^i$ is the union of the cells in $\mathcal C$ of dimension at most $i.$ We write $\stackrel{\circ}{X^i}$ for the union of cells in $\mathcal C$ of dimension $i.$ By convention, we write $X^{-1} = \emptyset.$ We say $\mathcal X$ is a stratification of $X$ if the closure of any $C\in\mathcal C$ in $X$ is a union of cells in $\mathcal C.$
\end{definition}
Note that the definition does not require or imply that $X$ is a closed subset of $\mathcal R^m$. This point has led to confusion. Unlike \cite{halup}, we restrict the definition to $C^2$ cell decompositions, since the existence of such stratifications is guaranteed by \cite{nguyen} and the arguments in \cite{halup} require refinement to the $C^2$ case anyways.
 The relations between $X, (X_i)_{i=0}^m,\stackrel{\circ}{X^i},$ and $\mathcal X$ in the definition will remain whenever we discuss stratifications.
Also, when a stratification is fixed, we write $T_xX$ for $T_x\stackrel{\circ}{X^i}$, where $x\in\stackrel{\circ}{X^i}.$ We write $P_x:\mathbb{R}^m\rightarrow T_xX$ for the orthogonal projection onto this space. The real inner product is fixed throughout.
\begin{definition}[\cite{halup}, Definition 1.2.3]
	Let $\mathcal X$ be a stratification and let $c,c', C,C'$ be positive elements of $\mathcal R$. We say a sequence $(a_i)_{i=0}^k$ with $a_i\in\stackrel{\circ}{X^{j_i}}$, $(j_i)_{i=0}^k$ a strictly decreasing sequence, is a plain $(c,c', C, C')-$chain for $\mathcal X$ if, for each $i$, $|a_0 - a_i| < cd(a_0,X^{j_i}),$ $d(a_0, X^{j_i-1})>Cd(a_0,X^{j_i})$, and, whenever $j<j_0$ and there is no $i$ such that $j_i = j$, $d(a_0, X^{j-1}) < c'd(a_0, X^j)$. An augmented $(c,c',C)-$chain is a sequence $(a_i)_{i=0}^k$ with $a_i\in\stackrel{\circ}{X^{j_i}}$, $(j_i)_{i=1}^k$ strictly decreasing, $j_0 = j_1$, and, for each $i>1$, $|a_0 - a_i| < cd(a_0,X^{j_i}),$ $d(a_0, X^{j_i-1})>Cd(a_0,X^{j_i})$, where we also have that, whenever $j<j_0$ and there is no $i$ such that $j_i = j$, $d(a_0, X^{j-1}) < c'd(a_0, X^j)$ and, finally, $d(a_0, a_1) < (C')^{-1}d(a_0,X^{j_0-1})$.  
\end{definition}
Note that increasing $c$ or $c'$ weakens the condition of being a $(c,c', C, C')-$chain, whereas increasing $C$ or $C'$ strengthens the condition. This is the meaning of the capitalization.
\begin{definition}[\cite{halup}, Definition 1.2.3]
	The stratification $\mathcal X$ is $(c,c', C,C',C'')-$Lipschitz if, whenever $(a_i)$ is a plain $(c,c',C,C')-$chain, $s\leq k$, we have $|(1-P_{a_0})P_{a_1}...P_{a_s}|\leq C''\frac{|a_0 - a_1|}{d(a_0,X^{j_s-1})}$, and whenever $(a_i)$ is an augmented $(c,c',C,C')-$chain, $|(P_{a_0}-P_{a_1})P_{a_2}...P_{a_s}|\leq C''\frac{|a_0 - a_1|}{d(a_0,X^{j_s-1})}$.  
	\\ \\
	If, for every $c>0$, there is a $C>0$ such that $\mathcal X$ is $(c,c,C,C,C)-$Lipschitz, we say $\mathcal X$ is a Lipschitz stratification.
\end{definition}
The proliferation of constants is useful. In earlier definitions of Lipschitz stratifications, the values of many of the constants were tied together, rendering the conditions non-monotonic; without this definition from \cite{halup}, we could not execute the proof of Lemma \ref{scale}.
\\ \\
We must now introduce the valuative Mostowski conditions. For the rest of this section, we fix $\mathcal R_0\prec \mathcal R$ and equip $\mathcal R$ with the valuation induced by the convex closure of $\mathcal R_0$ in $\mathcal R.$  We assume $\mathcal R$ is $|\mathcal R_0|^+-$saturated.
And we follow \cite{parus} in writing $v(d(a,\emptyset)) = -\infty.$ This abuse of notation shall cause no confusion.
\begin{definition}[\cite{halup}, Definition 1.6.1, Definition 1.6.8]
	Let $\mathcal X$ be a stratification of definable $X\subset\mathcal R^m$. A sequence $(a_i)_{i=0}^k$ is a plain val-chain if, with $a_i \in \stackrel{\circ}{X^{j_i}}$ for each $0<i\leq k$, $$v(|a_0 - a_i|) = v(d(a_0, X^{j_i})) = v(d(a_0,X^{j_{i-1}-1})) > v(d(a_0,X^{j_{i-1}})).$$
	It is an augmented val-chain if $j_0 = j_1$, the above inequality and equalities hold when $1<i\leq k$, and $v(|a_0-a_1|)<v(d(a_0, X^{j_0-1})).$ It is a weak plain val-chain if the equalities hold for all $0<i\leq k$, and it is a weak augmented val-chain if $j_0 = j_1$, the equalities hold for $1<i\leq k$, and $v(|a_0-a_1|)\geq v(d(a_0, X^{j_0-1})).$  
\end{definition}
The definitions of val-chains correspond to the definitions of $(c,c', C, C')-$chains; we now introduce the conditions on projections which correspond to the conditions given for $(c,c', C, C', C'')-$Lipschitz stratifications.
\begin{definition}
	Let $(a_i)_{i=0}^k$ be a sequence in $X$, $a_i\in\stackrel{\circ}{X^{j_i}}.$ Write $\lambda_i = v(|a_0 - a_i|)$ when $i\leq k$, and put $\lambda_{k+1} = v(d(a_0, X^{j_k-1})).$ If $(a_i)$ is a weak plain val-chain, the Mostowski condition for $(a_i)$ is that, for all $s\leq k$, $v(|(1-P_{a_0})P_{a_1}...P_{a_s}|)\geq \lambda_1 - \lambda_{s+1}.$
	If $(a_i)$ is a weak augmented val-chain, the Mostowski condition is that $v(|(P_{a_0}-P_{a_1})P_{a_2}...P_{a_s}|)\geq \lambda_1 - \lambda_{s+1}.$
	We say $\mathcal X$ is a valuative Lipschitz stratification if it satisfites the valuative Mostowski conditions at every val-chain and augmented val-chain.
\end{definition}
We are now ready to address Proposition 1.6.11 of \cite{halup}.
The proof of the implication (2')$\rightarrow$(3') given in \cite{halup} errs in its assumption that, when $a^0,...,a^n$ is a weak augmented val-chain, $a^1,...,a^n$ is a weak val-chain. Here is a corrected proof of their statement. 

\begin{lemma}\label{weak_chains}
Suppose $X\subset\mathcal{R}^{n}$ is definable with parameters. If $\mathcal X$ satisfies the valuative Mostowski conditions for every val-chain and augmented val-chain, then $\mathcal X$ satisfies the valuative Mostowski conditions for every weak val-chain and augmented weak val-chain.
\end{lemma} 
\begin{proof}
Let $a\in  \stackrel{\circ}{X^i}$. If $a = a_0, a_1, ..., a_k$ form a weak val-chain in $X$, removing those $a_i$ such that, if $j$ is such that $a_i = \stackrel{\circ}{X^j}$, then $v(d(a, X^{j-1})) = v(d(a, a_i)),$ with $d$ denoting the Euclidean distance function, yields a val-chain after re-indexing. Verifying the Mostowski conditions is then simple. 
	If $a = a_0, a_1,...,a_k$ is a weak augmented val-chain, we write $b = a^1$, remove such $a_i$ as above when $i>1$, and re-index. The resulting sequence  $a = a'_0, a'_1,...,a'_r$ is a plain val-chain. If $v(d(a,b))> v(d(a,X^{j-1}))$, then $a, b, a'_1,...,a'_r$ is an augmented val-chain and we are done. Assume that $v(d(a,b)) = v(d(a,X^{j-1})).$ 
\\ \\
Following the definition of a val-chain, we may find $b = b_0, b_1, ..., b_s$ such that $b_0, b_1,...,b_s,a'_2,...,a'_r$ is a plain val-chain. Such a sequence must also have $b_s, a'_1,a'_2,...,a'_r$ as an augmented val-chain. Write $\lambda_i = v(d(a, a_i'))$. Then our assumptions imply $v(|(P_{b_s} - P_{a_1'})P_{a_2'}...P_{a_r'}|) \geq v(d(b_s, a_1'))-\lambda_{r+1}\geq \lambda_1 - \lambda_{r+1}.$
We compute
$$P_{b_1}...P_{b_n}P_{a'_2}...P_{a'_r} - P_{b_s}P_{a_2}...P_{a'_r} = \sum_{i = 1}^{s-1} P_{b_i}...P_{b_s}P_{a'_2}...P_{a'_r} - P_{b_{i+1}}...P_{b_s}P_{a'_2}...P_{a'_r}$$
$$ = \sum_{i = 1}^{s-1} (P_{b_i}-1) P_{b_{i+1}}...P_{b_s}P_{a'_2}...P_{a'_r}$$
From the valuative Mostowski conditions, we have
$$ v( \sum_{i = 1}^{s-1} (P_{b_i}-1) P_{b_{i+1}}...P_{b_s}P_{a'_2}...P_{a'_r}) \geq \min_{1\leq i\leq s-1} \{\lambda_1 - \lambda_{r+1}\}$$
$$ =\lambda_1 - \lambda_{r+1}$$
By applying the Mostowski conditions to bound $|P_{b_1}...P_{b_s}P_{a'_2}...P_{a'_r}|$ and combining the three inequalities, we have $v((1-P_{b})P_{a'_1}...P_{a'_r})\geq \lambda_1 -\lambda_{r+1}$ which, when combined with the Mostowski conditions for $a, a'_1,...,a'_r$ yields the desired result.
\end{proof}
A corollary of Proposition 1.2.5 is
\begin{lemma}\label{scale}
If $\mathcal X = (X^i)_{i = 0}^m$ is a Lipschitz stratification of $X\subset\mathbb{R}^m$ definable with parameters, $U\subset\mathbb{R}^m$ is open and definable with parameters, $\mathcal Y$ is the stratification of $U\cap X$ induced by $\mathcal X$,  and there is some real constant $Q$ such that, for all $x\in U\cap X$ and $k\leq m$,  $d(x, Y^k) \leq Q d(x, X^k),$ then $\mathcal Y$ is a Lipschitz stratification.
\end{lemma}
\begin{proof}
	Applying criterion Proposition 1.2.5 (2) of Halupczok, we need to find an appropriate $C$ for every $c$. So let $c>1$ be arbitrary. By hypothesis, we can find a $C$ so that $\mathcal X$ satisfies the $(Qc, Qc, C, C, C)$ conditions. Now it is clear that $\mathcal Y$ satisfies the $(c, c, QC, QC, QC)$ Lipschitz conditions, as desired.
\end{proof}
We will need a way to obtain such a $U$ later. To this end, we demonstrate the following lemma, whose statement was inspired by communication with A. Parusinski:
\begin{lemma}\label{transverse}

	Suppose $\mathcal X = (X^i)_{i = 0}^m$ is a Lipschitz stratification of a closed $X\subset\mathbb{R}^m$ definable with parameters, and $U\subset\mathbb{R}^m$ is open and definable with parameters. Suppose also that $\partial U$ is a compact $C^2$ embedded submanifold of $\mathbb{R}^m$ with corners such that, for each $l\leq m$, the set $\partial_l \partial U$ of $l-$corners intersects each $\stackrel{\circ}{X^i}$ transversely. (The notation is unfortunate: $\partial U$ refers to the topological boundary of $U$, whereas $\partial_l$ refers to the codimension $l$ corner spaces, which have codimension $l+1$ in $\mathbb{R}^m$.) Write $\mathcal Y=(Y^k)_{k=0}^m$ for the stratification of $U\cap X$ induced by $\mathcal X$. Then there is some real constant $Q$ such that, for all $x\in U\cap X$ and $k\leq m$,  $d(x, Y^k) \leq Q d(x, X^k).$

\end{lemma}
\begin{proof}
	We obtain, for each $p\in \partial U$, a neighborhood $V$ of $p$ such that, whenever $x\in X\cap V$ and $y\in \partial U\cap V$, $T_xX + T_y\partial U = \mathbb{R}^m.$ (The argument for this is standard: if $p\notin X$, we take $V$ so $V\cap X=\emptyset$; if $p\in X$, let $\epsilon>0$ and take $C$ to be the cell in the cell decomposition inducing $\mathcal X$ such that $p\in C$ and write $n = \mathrm{dim}(C)$; our hypothesis gives $T_pX+T_p\partial U =\mathbb R^m$ and, since $p$ cannot be in the closure of any cell in the decomposition inducing $\mathcal X$ which has dimension less than or equal to $n$, we can take a neighborhood $V_0$ of $p$ disjoint from all such cells and we can find a neighborhood $V_1\subset V_0$ small enough that, whenever $q\in C\cap V$, $|P_p-P_q|<\epsilon$; we then find a neighborhood $V_2$ of $p$ such that, whenever $u\in V_2$, there is some $v\in V_1\cap C$ such that $|u - v| = d(u,X^{n})$, whence either condition (1) or (3) of Proposition 1.6.5 in \cite{halup} give constants $C''$ for certain $(\alpha,\beta,\gamma,\delta)$ where every $u\in V_2$ is the first point in a $(\alpha,\beta,\gamma,\delta)-$chain in which an element of $C\cap V_1$ appears such that $\mathcal X$ is a $(\alpha,\beta,\gamma,\delta,C'')-$Lipschitz stratification, so if the diameter of $V_3\subset V_2$ is less than $\frac{\epsilon }{2kC''}$ and $2(V_3 - p)\subset V_2 - p$, the set operations having their usual meaning, the telescoping sum analysis of the projections we used in the proof of Lemma \ref{weak_chains} gives, for each $z\in V_3\cap X$, $|(1-P_z)P_p|<\epsilon$; likewise, we may pick a neighborhood $W_0$ of $p$ such that there is a definable $C^2$ map $\phi$ from $W_0\cap \partial U$ onto an open subset of $(\mathbb R_{\geq 0})^s$, where we shall write coordinates $x^1,...,x^s$, for some $s$ which has a $C^2$ inverse, and indeed $\phi^{-1}$ being $C^1$ is enough to guarantee there is some neighborhood $\mathcal O$ of $\phi(p)$ in $(\mathbb R_{\geq 0})^s$ on which, whenever the $\nu-$th coordinate of $\phi(p)$ is nonzero and $t\in\mathcal O$, $|\frac{\partial\phi^{-1}}{\partial x^\nu}|_t- \frac{\partial\phi^{-1}}{\partial x^\nu}|_{\phi(p)}|<\epsilon$; put $W_1 = \phi^{-1}(\mathcal O)$ and note that, as $\epsilon$ was arbitrary, the result follows.)
	In fact, the parenthetical argument gives us more: from it and the compactness of $\partial U$, we obtain some $r,\rho>0$ such that, whenever $q\in\partial U, x\in X,$ and $|q - x|<r$, we have $\inf_{|v| = 1}[|P_qv|+|P_xv|] > \rho$.
	We write $V^\perp$ for the orthogonal complement of a vector subspace $V\leq \mathbb{R}^m.$
	Put $r' = \rho r.$ We claim that, whenever $x\in X^k$ and $d(x,\partial U) < r'$, there is $y\in X^k\cap \partial U$ such that $|x - y |<\rho^{-1} d(x,\partial U).$ The statement of the lemma follows immediately from this claim.
	\\ \\
	Let $x_0\in X^k$ with $d(x_0,\partial U) < r'$ be arbitrary, and take $y_0\in\partial U$ so $|x_0 - y_0 | = d(x_0,\partial U).$ Suppose $y_0\in\partial_l\partial U.$ We may suppose by induction that the claim holds for all such $(x,y)$ satisfying these conditions with $x\in X^{k'}, y \in \partial_{l'}\partial U$ when $k'\leq k, l'\geq l,$ and $k'<k$ or $l'>l.$ Then define a vector field $Z$ on $\partial_l\partial U\times X^k\backslash [(\partial U \cap X^k)\times(\partial U \cap X^k)]$ by $Z_{(s,t)} = (P_s\frac{t-s}{|t-s|},P_t\frac{s-t}{|s-t|}).$ Writing $d:\partial_l\partial U \times X^k\rightarrow \mathbb{R}$ for the distance function, whenever $(s,t)$ is in the domain of $Z$, the inequality $Z_{(s,t)}d\leq -\rho$ holds. Let $\gamma$ be a maximal integral curve for $Z$ with $\gamma(0) = (y_0,x_0)$. Suppose $\gamma$ is defined on the interval $(a,b).$ (Of course $b<\infty.$) Then $\lim_{t\rightarrow b}\gamma(t)$ is either in $\partial_{l+1}\partial U\times X^k$, $\partial_l\partial U\times X^{k-1}$, or $(\partial U\cap X^k)\times (\partial U \cap X^k).$ In the former cases, the inductive hypothesis yields the result, while in the latter, the conclusion is obvious.
\end{proof}
\section{Isotopy and Finiteness of Volume}
We give here a result of Nguyen. A complete proof can be found in \cite{nguyen}, though the statement Nguyen gives is weaker.
\begin{theorem} (\cite{nguyen}, Theorem 4.4.5)\label{ng_strat}
	If $X\subset\mathcal R^m$ is closed and definable with parameters, and $(A_i)_{i<k}$ is a finite sequence of definable subsets of $X,$ then there is a Lipschitz stratification $\mathcal X = (X^i)_0^m$ of $X$ such that every $A_i$ is a finite union of connected components of $\stackrel{\circ}{X^i}$ for various $i.$
\end{theorem}
We are now in a position to prove our main result. The proof is a combination of the proof of Theorem 1.6 and Lemma 1.7 of \cite{parus} and Proposition 1 of \cite{clr}.
\begin{theorem}
	$L(f)$ is definable with parameters in $\mathcal R.$
\end{theorem}
\begin{proof}
	As indicated in the introduction, we simplify the notation by determining which fibers $A_x$ of $A\subset\mathbb R^{m+n}$ have finite $n-$volume. Without loss of generality, assume each $A_x$ is closed.
Let $\alpha:\mathbb R^{n}\backslash \{0\}\rightarrow\mathbb R^{n} \backslash \{0\}$ denote the inversion about the unit sphere: $\alpha(v) = |v|^{-2}v.$ We may assume that the intersection of each $A_x$ with the unit ball in $\mathbb R^{n}$ is empty. Take $B = \{(x,y)\in\mathbb R^{m+n}:(x,\alpha(y))\in A\}\cup\mathbb R^m\times\{0\}.$ By Theorem \ref{ng_strat}, we may take a definable Lipschitz stratification $\mathcal X$ of $B$ compatible with $\mathbb R^m \times\{0\}.$ 
	The set $S$ of $x$ such that $x$ is a regular value of $\pi_m\upharpoonright \stackrel{\circ}{X^k}$ for each $k$, where $\pi_m$ denotes projection onto the $\mathbb R^m$ coordinate, is open and definable. By induction on the dimension $m$, we need only concern ourselves with $A_x$ for $x\in S.$ About each $x\in S$,there is a convex neighborhood $V$ with $\overline V\subset S;$ fix such $x_0$ and $V.$ By Lemmas \ref{scale} and \ref{transverse}, $\mathcal X\cap \pi_m^{-1}(V)$ is a Lipschitz stratification of $B\cap \pi_m^{-1}(V).$ 
	\\ \\
	We now follow the proof of Lemma 1.7 in \cite{parus} closely. We claim there is a neighborhood of $x_0$ such that every Lipschitz vector field on that neighborhood can be lifted to a Lipschitz vector field on its pre-image under $\pi_m.$ It suffices to check this for every vector field in a $C^1$ frame on a neighborhood of $x_0,$ as the coefficients of any Lipschitz vector field with respect to such a frame are Lipschitz functions and a linear combinations of lifts is a lift of linear combinations. For simplicity, then, let us lift a coordinate vector field $\partial_i.$ By the compactness of $\pi_m^{-1}(x_0)$ and the choice of $V$, we know there is some constant $L$ such that, for all $v\in\mathbb R^m$ and $\tilde x\in\pi_m^{-1}(x)$, there is some $w\in\mathbb R^{m+n}$ such that $w\in T_{\tilde x} B,\pi_m(w) = v,$ and $|w| <L.$ For any $\tilde x\in\pi_m^{-1}(x)$, we can find a neighborhood $N$ of $\tilde x$ in the cell $C$ of $\mathcal X$ containing $\tilde x$ on which some Lipschitz vector field $Y$ with Lipschitz constant $G$ has $|Y_p|<L$ at every $p$ where $Y$ is defined and ${\pi_m}_*Y = \partial_i.$ Take $\psi$ to be a $C^1$ definable bump function which is identically $1$ on a neighborhood of $\tilde x$ and which is supported in $N.$ Then $\psi Y$ has Lipschitz constant $G+\sup |\nabla \psi|$ on $N$, and extending $\psi$ to $X^{\mathrm{dim}(C)}$ by declaring it to be 0 elsewhere yields a Lipschitz vector field on that skeleton. Now Proposition 1.3 of \cite{parus} guarantees that this can be extended to a Lipschitz vector field on $X.$ Having done this for each $\partial_i$, obtaining a corresponding $\xi_i$, we note that the matrix $[{\pi_m}_*\xi_1|...|{\pi_m}_*\xi_m]$ is a Lipschitz function which, at $\tilde x$, equals the identity matrix. Therefore, the inverse matrix $M$ is also a Lipschitz function on a neighborhood of $\tilde x$ in $X$, and $\xi_iM^i_j$ is a lift of $\partial_i.$ By covering $\pi_m^{-1}(x)$ with the neighborhoods resulting from applying this method to various $\tilde x$'s and taking a finite subcovering, we may take a definable partition of unity subordinate to such a covering and take the corresponding weighted sum to obtain a lift of $\partial_i$ to a neighborhood of $\pi_m^{-1}(x).$  This proves the claim.
	\\ \\
	Take a convex neighborhood $W$ of $x$ such that every Lipschitz vector field on $W$ lifts to a Lipschitz vector field on $\pi_m^{-1}(x).$ Let $y_0\in W\backslash \{x_0\}$ be arbitrary and let $Z$ be the vector field on $V\backslash \{x_0\}$ given by $Z_z = \frac{x_0-y_0}{|x_0-y_0|}$. Take $\tilde{Z}$ on $\pi^{-1}(V\backslash\{x\})\cap B$ which is Lipschitz and compatible with $\mathcal X$ (meaning if $s\in\stackrel{\circ}{X^k}$, then $\tilde{Z}_s\in T_s\stackrel{\circ}{X^k}$) such that ${\pi_m}_*\tilde{Z} = Z$. The flow of $\tilde{Z}$ gives a Lipschitz homeomorphism between $B_{y_0}$ and $B_{x_0}$. Reversing the roles of $x_0$ and $y_0$ and considering the flow of $-\tilde{Z}$ gives the Lipschitz homeomorphism from $B_{x_0}$ to $B_{y_0}$. So $B_{x_0}$ and $B_{y_0}$ are bi-Lipschitz homeomorphic and, because $\mathcal X$ is compatible with $\mathbb R^m\times\{0\}$, the flows and hence the homeomorphisms fix $0\in\mathbb R^{n}.$ Now Affirmation 2 of \cite{clr} yields that $B_{x_0}$ has finite volume if and only if $B_{y_0}$ does. Therefore, $L(f)\cap S$ is a union of connected components of $S$, and the result follows.

\end{proof}
\bibliography{polybib}
\bibliographystyle{ieeetr}
\end{document}